\documentclass[a4paper,11pt]{article}

\usepackage{amsfonts}
\usepackage{amsmath}
\usepackage{amssymb}
\usepackage{graphicx}
\usepackage{amscd}
\usepackage{tikz}
\usepackage{epsfig}
\usepackage{psfrag}

\title{A new solution to square matrix completion problem}

\author{Marija Dodig\thanks{CEAFEL, Departamento de Mat\'ematica, Universidade de Lisboa, Edificio C6, Campo Grande, 
1749-016 Lisbon, Portugal, and
Mathematical Institute SANU, Knez Mihajlova 36, 11000 Belgrade, Serbia. ({\tt msdodig@fc.ul.pt}). Corresponding author.}
\and
Marko Sto\v si\'c\thanks{CAMGSD, Departamento de Matem\'atica, 
Instituto Superior T\'ecnico,
Av. Rovisco Pais 1, 1049-001 Lisbon, Portugal, and
Mathematical Institute SANU, Knez Mihajlova 36, 11000 Belgrade, Serbia. }
}

\date{}

\newtheorem{theorem}{Theorem}

\newtheorem{definition}{Definition}
\newtheorem{example}{Example}

\newtheorem{lemma}{\indent Lemma}

\newtheorem{problem}{Problem}

\newtheorem{remark}{Remark}

\def\max{\mathop{\rm max}}

\def\lcm{\mathop{\rm lcm}}

\def\kraj{\hfill\rule{6pt}{6pt}}

\def\deg{\mathop{\rm deg}}

\def\F{\mathbb{F}}

\begin{document}

\maketitle

\begin{abstract}
In this paper we give a novel solution to a classical completion problem for square matrices. This problem was studied by many authors through time, and it is completely solved in \cite{super,pol}. 
In this paper we relate this classical problem to a purely combinatorial question involving partitions of integers and their majorizations studied in \cite{silvaL}. We show surprising relations in these approaches and  as a corollary, we obtain a new combinatorial result on partitions of integers.
\end{abstract}



\section{Introduction}

In this paper we consider the following classical matrix completion problem:
\begin{problem}\label{p}
Describe the possible similarity class of a square matrix with a prescribed submatrix.\end{problem}
Problem \ref{p} has a long history - it is one of the most studied matrix completion problems. Various particular cases have beed solved, see e.g. \cite{oliveira,sa,tomp,zabala}. The necessary conditions for it were obtained by Gohberg,  Kaashoek, and van Schagen in \cite{face}.  Significantly more difficult is proving the sufficiency of the conditions from \cite{face}. First attempt of proving sufficiency was made by Cabral  and Silva in \cite{cabral}, where an implicit solution to Problem \ref{p} was obtained. Later on in \cite{super}, Dodig and Sto\v si\'c  gave a complete, explicit and constructive solution to Problem \ref{p} \cite[Theorem 1]{super}. Recently, in \cite{pol} a new, purely combinatorial  and more direct and elegant way to solve Problem \ref{p} was given in \cite[Corollary 5]{pol}. In fact, in \cite[Section 4]{pol} (see also \cite{super}) has been shown that
Problem \ref{p} has a solution if and only if the following theorem is valid. Throughout the paper $\F$ is an algebraically closed field.
\begin{theorem}\label{T}\cite{super,pol} Let  $\tilde{\alpha}:\tilde{\alpha}_1|\cdots|\tilde{\alpha}_n$ and $\tilde{\gamma}:\tilde{\gamma}_1|\cdots|\tilde{\gamma}_{n+m+p}$ 
be chains of homogeneous polynomials  from $\F[\lambda,\mu]$, and let $c_1\ge\cdots\ge c_m$ and $r_1\ge\cdots\ge r_p$ be nonnegative integers, such that 
\begin{eqnarray*}
(i)&\tilde{\gamma}_i\mid\tilde{\alpha}_i\mid\tilde{\gamma}_{i+m+p},\quad i=1,\ldots,n,\label{(i1)}\\
(ii)&(c_1+1,\ldots,c_m+1)\cup(r_1+1,\ldots,r_p+1)\prec(d(\tilde{\sigma}_{m+p}(\tilde{\alpha},\tilde{\gamma})),\ldots, d(\tilde{\sigma}_{1}(\tilde{\alpha},\tilde{\gamma}))).\label{(ii1)}
\end{eqnarray*} 
Then there exists a chain of homogeneous polynomials
$\tilde{\beta}:\tilde{\beta}_1|\cdots|\tilde{\beta}_{n+m}$, from $\F[\lambda,\mu]$ which satisfies:\begin{eqnarray}
&\tilde{\beta}_i|\tilde{\alpha}_i|\tilde{\beta}_{i+m},\quad i=1,\ldots,n,\label{p1p}\\
&\tilde{\gamma}_i|\tilde{\beta}_i|\tilde{\gamma}_{i+p},\quad i=1,\ldots,n+m,\label{p2}\\
&(c_1+1,\ldots,c_m+1)\prec(d(\tilde{\sigma}_m(\tilde{\alpha},\tilde{\beta})),\ldots,d(\tilde{\sigma}_1(\tilde{\alpha},\tilde{\beta}))),\\
&(r_1+1,\ldots,r_p+1)\prec(d(\tilde{\sigma}_p(\tilde{\beta},\tilde{\gamma})),\ldots,d(\tilde{\sigma}_1(\tilde{\beta},\tilde{\gamma}))).\label{n4p}
\end{eqnarray}




\end{theorem}

\vspace{0,5cm}

\noindent Here for any two polynomial chains $\tilde{\delta}: \tilde{\delta}_1|\cdots|\tilde{\delta}_x$ and $\tilde{\epsilon}:\tilde{\epsilon}_1|\cdots|\tilde{\epsilon}_{x+y}$ such that $\tilde{\epsilon}_i|\tilde{\delta}_i|\tilde{\epsilon}_{i+y}$, $i=1,\ldots,x$, we define:
$$\tilde{\sigma}_i(\tilde{\delta},\tilde{\epsilon})={{\tilde{\pi}_i}(\tilde{\delta},\tilde{\epsilon})\over{\tilde{\pi}_{i-1}}(\tilde{\delta},\tilde{\epsilon})},\quad \tilde{\pi}_i(\tilde{\delta},\tilde{\epsilon})=\prod_{j=1}^{x+i}{\lcm(\tilde{\delta}_{j-i}, \tilde{\epsilon}_j)}, \quad i=0,\ldots,y.$$


Thus, in order to solve Problem \ref{p}, we are left with proving Theorem \ref{T}.  This has been done in two completely different ways in \cite{super} and in \cite{pol}. However, in this paper we present
another solution to Problem \ref{p}. We study and show surprising equivalence between Theorem \ref{T}
 and combinatorial results on majorization of partitions  obtained in \cite{silvaL}. 
 As a corollary of these relations  we obtain a new combinatorial result on majorization of partitions in Lemma \ref{nova}.

 \section{Notation}
  Throughout the paper we deal with (chains of) homogeneous polynomials
from $\F[\lambda,\mu]$. By homogeneous irreducible factors of the homogeneous polynomial $\tilde{f}\in\F[\lambda,\mu]$, we mean homogeneous irreducible polynomials from $\F[\lambda,\mu]$ that divide $\tilde{f}$. For  a polynomial chain $\tilde{\alpha}_1|\cdots|\tilde{\alpha}_n$, we assume $\tilde{\alpha}_i\in\F[\lambda,\mu]$  are all monic, nonzero polynomials. By convention we have $\tilde{\alpha}_i=1$, for $i\le 0$, and $\tilde{\alpha}_i=0$, for $i\ge n+1$.
 Also, we assume 
$d(1)=\deg(1)=0$, and $d(0)=+\infty$.

 By a partition of integers, we mean a non-increasing sequence of non-negative integers. For a partition 
$\mathbf{a}=(a_1,\ldots,a_m)$, we assume $a_{m+1}=a_{m+2}=\cdots=0$, and we identify two partitions differing only by a tail of zeros. Also, by $|\mathbf{a}|$ we denote $\sum_{i=1}^m{a_i}$, 
{{and by $\overline{\mathbf{a}}=(\overline{a}_1,\ldots,\overline{a}_{{|\mathbf{a}|}})$ we denote the dual partition of $\mathbf{a}$. Here $\overline{a}_i=\sharp\{j|a_j\ge i\}$, $i=1,\ldots,{{|\mathbf{a}|}}$.}} 

For any two partitions $\mathbf{a}=(a_1, a_2,\ldots)$, and $\mathbf{b}=(b_1,b_2,\ldots)$ with $a_i\ge b_i$, $i\ge 1,$ by $\mathbf{a}-\mathbf{b}$ we denote a partition obtained by ordering the elements $a_i-b_i$, $i\ge1$, in the non-increasing order. Also, we put $\mathbf{a}+\mathbf{b}=(a_1+b_1,a_2+b_2,\ldots)$. The partition $\mathbf{a}\cup\mathbf{b}$ is defined as a partition whose non-zero elements are precisely the non-zero elements of partitions $\mathbf{a}$ and $\mathbf{b}$ ordered in non-increasing order. Recall that 
$$\overline{\mathbf{a}\cup\mathbf{b}}=\overline{\mathbf{a}}+\overline{\mathbf{b}}.$$ 
We also recall the definition of the classical majorization, \cite{ped}: 
 \begin{definition}Let {
$\mathbf{a}=(a_1,a_2,\ldots,a_n)$} and  
{$\mathbf{b}=(b_1,b_2,\ldots,b_n)$}
be two sequences of nonnegative integers, not necessarily non-increasing. Let $\sigma^1$ and $\sigma^2$ be two permutations of the set $\{1,\ldots,n\}$ such that 
$a_{\sigma^1(1)}\ge a_{\sigma^1(2)}\ge  \cdots\ge a_{\sigma^1(n)}$
and
$b_{\sigma^2(1)}\ge b_{\sigma^2(2)}\ge \cdots b_{\sigma^2(n)}.$

If
$$\sum_{i=1}^n{a_i} = \sum_{i=1}^n{b_i},$$ 
and
$$\sum_{i=1}^j{a_{\sigma^1(i)}}\le\sum_{i=1}^j{b_{\sigma^2(i)}},\quad j=1,\dots,n-1,$$
then we say that $\mathbf{a}$ is {\it{majorized}} by $\mathbf{b}$, and write {$\mathbf{a}\prec\mathbf{b}$}.

\end{definition}
We note  that $\mathbf{a}\prec\mathbf{b}$ is equivalent to $\overline{\mathbf{b}}\prec\overline{\mathbf{a}}$, and also if $\mathbf{a}\prec\mathbf{b}$ and $\mathbf{b}\prec\mathbf{c}$, then $\mathbf{a}\prec\mathbf{c}$.


\section{Combinatorial lemmas}

 In \cite{silvaL} we have studied series connections of arbitrarily many linear systems. 
  As the main result, we completely determined the controllability and the possible controllability indices of a system obtained by a special series connection of arbitrarily many linear systems. As the crucial part of the proof of the main result in \cite{silvaL},  we have obtained the following combinatorial result involving classical majorizations of partitions of integers:
 
  \begin{lemma}\cite[Lemma 5]{silvaL}\label{silvaL}
  Let $d_1\ge\cdots\ge d_s\ge 0$ and $t_1\ge\cdots\ge t_s\ge 0$ be nonincreasing sequences of nonnegative integers, such that  $d_i\ge t_i$, $i=1,\ldots,s$. Let $A_1\ge\cdots\ge A_s\ge 0$ and $B_1\ge\cdots\ge B_s\ge 0$ be nonincreasing sequences of nonnegative integers such that 
 \begin{equation}(d_1-t_1,\ldots,d_s-t_s)\prec (A_1+B_1,\ldots,A_s+B_s).\end{equation}
 Then there exists a nonincreasing sequence $f_1\ge\cdots\ge f_s$ of non-negative integers such that 
 \begin{equation} d_i\ge f_i\ge t_i,\quad i=1,\ldots,s,\end{equation}
 and such that 
 \begin{equation}(f_1-t_1,\ldots,f_s-t_s)\prec (A_1,\ldots,A_s),\end{equation}
 \begin{equation}(d_1-f_1,\ldots,d_s-f_s)\prec (B_1,\ldots,B_s).\end{equation}
  \end{lemma}
\begin{remark}In the original formulation of this lemma in \cite{silvaL}, it was required that  $d_i$, $t_i$ and $f_i$ are strictly positive for all $i=1,\ldots,s$, --- this was motivated by the particular completion problem that this was related to. However, it is clear that the conditions of the lemma depend only on the differences $d_i-t_i$ and therefore clearly remain valid if one increases (or decreases) all $d_i$'s, $f_i$'s, and $t_i$'s by the same value. Therefore one can assume that all $d_i$'s, $f_i$'s, and $t_i$'s are nonnegative integers.\end{remark}

In this paper we shall show remarkable relationship between Lemma \ref{silvaL} and Theorem \ref{T}, see Remark \ref{r}.  Moreover, inspired by this relation we  give a new combinatorial result on partitions of integers and their majorizations that we show to be equivalent to Theorem \ref{T}. It is a very surprising connection between two completely unrelated problems. This novel, general and interesting combinatorial result is given in the following lemma:

\begin{lemma}\label{nova}
 Let $\mathbf{d}^i=(d^i_1,\ldots,d^i_s)$ and $\mathbf{t}^i=(t^i_1,\ldots,t^i_s)$, $i=1,\ldots,k$, be partitions of nonnegative integers, such that
 $d^i_j\ge t^i_j$, $j=1,\ldots,s$, $i=1,\ldots,k$.  Let 
$\mathbf{A}=(A_1,\ldots,A_s)$ and $\mathbf{B}=(B_1,\ldots,B_s)$ be partitions of  nonnegative integers
 such that
 \begin{equation}(\mathbf{d}^1-\mathbf{t}^1)\cup\cdots\cup(\mathbf{d}^k-\mathbf{t}^k)\prec \mathbf{A+B}.
 \end{equation}
 Then there exist partitions $\mathbf{f}^i=(f^i_1,\ldots, f^i_s)$, $i=1,\ldots,k,$ of nonnegative integers such that 
\begin{eqnarray}
&&d^i_j\ge f^i_j\ge t^i_j,\quad \textrm{ for all } \quad i=1,\ldots,k,\quad j\ge 1\label{5}\\
&&(\mathbf{f}^1-\mathbf{t}^1)\cup(\mathbf{f}^2-\mathbf{t}^2)\cup\cdots\cup(\mathbf{f}^k-\mathbf{t}^k)\prec \mathbf{A}\\
&&(\mathbf{d}^1-\mathbf{f}^1)\cup(\mathbf{d}^2-\mathbf{f}^2)\cup\cdots\cup(\mathbf{d}^k-\mathbf{f}^k)\prec \mathbf{B}.\label{7}
\end{eqnarray}
\end{lemma}

Clearly, for $k=1$ Lemma \ref{nova} reduces to Lemma \ref{silvaL}.

\section{A new proof of Problem \ref{p}}

Before proceeding with our main result, let us introduce some notation.  Let  $\tilde{\alpha}:\tilde{\alpha}_1|\cdots|\tilde{\alpha}_n$ and $\tilde{\gamma}:\tilde{\gamma}_1|\cdots|\tilde{\gamma}_{n+m+p}$ 
be polynomial chains of homogeneous polynomials from $\F{[\lambda,\mu]}$, and let $c_1\ge\cdots\ge c_m$ and $r_1\ge\cdots\ge r_p$ be nonnegative integers. Let $\psi_1,\ldots,\psi_k$ be irreducible  factors of $\tilde{\gamma}_{n+m+p}$.  For every $i=1,\ldots,k$, let $\mathbf{a}^i=(a^i_1,\ldots,a^i_n)$ and $\mathbf{g}^i=(g^i_1,\ldots,g^i_{n+m+p})$ be partitions corresponding to the  $\psi_i$ elementary divisor of the polynomial chains  $\tilde{\alpha}:\tilde{\alpha}_1|\cdots|\tilde{\alpha}_n$ and $\tilde{\gamma}:\tilde{\gamma}_1|\cdots|\tilde{\gamma}_{n+m+p}$, respectively.
More precisely:
 $$\tilde{\alpha}_i=\psi_1^{a^1_{n+1-i}}\psi_2^{a^2_{n+1-i}}\ldots\psi_k^{a^k_{n+1-i}},\quad i=1,\ldots,n,$$
 $$\tilde{\gamma}_i=\psi_1^{g^1_{n+m+p+1-i}}\psi_2^{g^2_{n+m+p+1-i}}\ldots\psi_k^{g^k_{n+m+p+1-i}},\quad i=1,\ldots,n+m+p.$$

 Then, if $$\tilde{\gamma}_i|\tilde{\alpha}_i|\tilde{\gamma}_{i+m+p},\quad i=1,\ldots,n,$$
from the definition of $\tilde{\sigma}(\tilde{\alpha},\tilde{\gamma}), \mathbf{g}^1, \ldots, \mathbf{g}^k,\mathbf{a}^1$,\ldots, $\mathbf{a}^k$, we have
\begin{equation}\label{1}(d(\tilde{\sigma}_{m+p}(\tilde{\alpha},\tilde{\gamma})),\ldots, d(\tilde{\sigma}_{1}(\tilde{\alpha},\tilde{\gamma})))=d(\psi_1)\overline{\overline{\mathbf{g}^1}-\overline{\mathbf{a}^1}}+\cdots+d(\psi_k)\overline{\overline{\mathbf{g}^k}-\overline{\mathbf{a}^k}}.\end{equation}

Since $\F$ is algebraically closed field, we have that $d(\psi_i)=1$, $i=1,\ldots,k$, i.e. (\ref{1}) is equal to
\begin{equation}\label{1u}(d(\tilde{\sigma}_{m+p}(\tilde{\alpha},\tilde{\gamma})),\ldots, d(\tilde{\sigma}_{1}(\tilde{\alpha},\tilde{\gamma})))=\overline{\overline{\mathbf{g}^1}-\overline{\mathbf{a}^1}}+\cdots+\overline{\overline{\mathbf{g}^k}-\overline{\mathbf{a}^k}}.\end{equation}

Now we can give our main result:

\begin{theorem} \label{main}
Theorem \ref{T} is equivalent to Lemma \ref{nova}.
\end{theorem}
\noindent\textbf{Proof:} 
We  start by proving that Lemma \ref{nova} implies Theorem \ref{T}. 

Thus, let the conditions $(i)$ and $(ii)$ from Theorem \ref{T}
be valid. Then $(ii)$ and (\ref{1u}) together give
\begin{equation}\label{2uun}
{\overline{\overline{\overline{\mathbf{g}^1}-\overline{\mathbf{a}^1}}+\cdots+\overline{\overline{\mathbf{g}^k}-\overline{\mathbf{a}^k}}}}\prec\overline{(c_1+1,\ldots,c_m+1)\cup(r_1+1,\ldots,r_p+1)}\end{equation}

Let $\mathbf{c}=(c_1+1,\ldots,c_m+1)$ and let $\mathbf{r}=(r_1+1,\ldots,r_p+1)$, and let $\mathbf{A}:=\overline{\mathbf{c}}$ and $\mathbf{B}:=\overline{\mathbf{r}}$, $\mathbf{A}=(A_1,A_2,\ldots)$ and $\mathbf{B}=(B_1,B_2,\ldots)$. Then $A_1=m$ and $B_1=p$. Since $\overline{\mathbf{c}\cup\mathbf{r}}=\overline{\mathbf{c}}+\overline{\mathbf{r}}$, we have that 
$$\overline{(c_1+1,\ldots,c_m+1)\cup(r_1+1,\ldots,r_p+1)}=\mathbf{A}+\mathbf{B}.$$
 Thus, (\ref{2uun}) is equal to
\begin{equation}\label{2n}
\overline{\overline{\overline{\mathbf{g}^1}-\overline{\mathbf{a}^1}}+\cdots+\overline{\overline{\mathbf{g}^k}-\overline{\mathbf{a}^k}}}\prec \mathbf{A}+\mathbf{B}.
\end{equation}
\noindent Let us denote
 by $\mathbf{d}^i:=\overline{\mathbf{g}^i}$, with $\mathbf{d}^i=(d^i_1,d^i_2,\ldots)$, and $\mathbf{t}^i:=\overline{\mathbf{a}^i}$ with $\mathbf{t}^i=(t^i_1,t^i_2,\ldots)$, $i=1,\ldots,k.$ Then (\ref{2n}) becomes
\begin{equation}\label{4n}
(\mathbf{d}^1-\mathbf{t}^1)\cup(\mathbf{d}^2-\mathbf{t}^2)\cup\cdots\cup(\mathbf{d}^k-\mathbf{t}^k)\prec \mathbf{A}+\mathbf{B}.\end{equation}
By Lemma \ref{nova}, there exist  partitions $\mathbf{f}^1,\ldots,\mathbf{f}^k$ such that
\begin{eqnarray}
&&d^i_j\ge f^i_j\ge t^i_j,\quad \textrm{ for all } \quad i=1,\ldots,k,\quad j\ge 1\label{51}\\
&&(\mathbf{f}^1-\mathbf{t}^1)\cup(\mathbf{f}^2-\mathbf{t}^2)\cup\cdots\cup(\mathbf{f}^k-\mathbf{t}^k)\prec \mathbf{{A}}\label{61}\\
&&(\mathbf{d}^1-\mathbf{f}^1)\cup(\mathbf{d}^2-\mathbf{f}^2)\cup\cdots\cup(\mathbf{d}^k-\mathbf{f}^k)\prec \mathbf{B}.\label{71}
\end{eqnarray}
Let  $$\mathbf{b}^i=\overline{\mathbf{f}^i},\quad i=1,\ldots,k,$$ and let $\tilde{\beta}:\tilde{\beta}_1|\cdots|\tilde{\beta}_{n+m}$, be a polynomial chain such that the only irreducible factors of $\tilde{\beta}_{n+m}$ are $\psi_1,\ldots,\psi_k$, and such that for all $i=1,\ldots,k$, $\mathbf{b}^i$ is the partition corresponding to the $\psi_i$ elementary divisor of $\tilde{\beta}$, i.e.
$$\tilde{\beta}_i=\psi_1^{b^1_{n+m+1-i}}\psi_2^{b^2_{n+m+1-i}}\ldots\psi_k^{b^k_{n+m+1-i}},\quad i=1,\ldots,n+m.$$
From (\ref{61}) we have that for every $i=1,\ldots,k,$ it is valid that $\max\{f^i_j-t^i_j | j\ge 1\}\le A_1=m$, and so $\overline{t^i}_j \ge \overline{f^i}_{j+m}$, for all $i$ and $j$.
Analogously from (\ref{71}) we have that for every $i=1,\ldots,k,$ it is valid that $\max\{d^i_j-f^i_j | j\ge 1\}\le B_1=p$, and so $\overline{f^i}_j \ge \overline{d^i}_{j+p}$, for all $i$ and $j$. Together with (\ref{51}) this gives
$$b^i_j\ge a^i_j\ge b^i_{j+m}, \quad \textrm{ and }\quad g^i_j\ge b^i_j\ge g^{i}_{j+p}, \quad i=1,\ldots,k,\quad  j\ge 1.$$
Hence,
\begin{eqnarray}
&\tilde{\beta}_i|\tilde{\alpha}_i|\tilde{\beta}_{i+m},\quad i=1,\ldots,n,\label{8}\\
&\tilde{\gamma}_i|\tilde{\beta}_i|\tilde{\gamma}_{i+p},\quad i=1,\ldots,n+m.\label{9}\end{eqnarray}

 Then the duals of (\ref{61}) and (\ref{71})  give
$$(c_1+1,\ldots,c_m+1)=\overline{\mathbf{A}}\prec \overline{(\overline{\mathbf{b}^1}-\overline{\mathbf{a}^1})\cup\cdots\cup(\overline{\mathbf{b}^k}-\overline{\mathbf{a}^k})}=$$
$$=\overline{(\overline{\mathbf{b}^1}-\overline{\mathbf{a}^1})}+\cdots+\overline{(\overline{\mathbf{b}^k}-\overline{\mathbf{a}^k})}=(d(\tilde{\sigma}_m(\tilde{\alpha},\tilde{\beta})),\ldots,d(\tilde{\sigma}_1(\tilde{\alpha},\tilde{\beta})))$$
and 
$$(r_1+1,\ldots,r_p+1)=\overline{\mathbf{B}}\prec \overline{(\overline{\mathbf{g}^1}-\overline{\mathbf{b}^1})\cup\cdots\cup(\overline{\mathbf{g}^k}-\overline{\mathbf{b}^k})}=$$
$$=\overline{(\overline{\mathbf{g}^1}-\overline{\mathbf{b}^1})}+\cdots+\overline{(\overline{\mathbf{g}^k}-\overline{\mathbf{b}^k})}=(d(\tilde{\sigma}_p(\tilde{\beta},\tilde{\gamma})),\ldots,d(\tilde{\sigma}_1(\tilde{\beta},\tilde{\gamma}))).$$
Hence, such defined $\tilde{\beta}_i$'s  satisfy (\ref{p1p})--(\ref{n4p}), as wanted.\\

Now, suppose that Theorem \ref{T} is valid, and let us prove Lemma \ref{nova}. 
Let $\mathbf{d}^i=(d^i_1,\ldots,d^i_s)$, $\mathbf{t}^i=(t^i_1,\ldots,t^i_s)$, $i=1,\ldots,k$, and let $\mathbf{A}=(A_1,\ldots,A_s)$ and $\mathbf{B}=(B_1,\ldots,B_s)$ be partitions such that
\begin{equation}{d}^i_j\ge {t}^i_j,\quad i=1,\ldots,k,\quad  j=1,\ldots,s,\label{?1}\end{equation} and such that 
\begin{equation}\label{*1}
(\mathbf{d}^1-\mathbf{t}^1)\cup(\mathbf{d}^2-\mathbf{t}^2)\cup\cdots\cup(\mathbf{d}^k-\mathbf{t}^k)\prec \mathbf{A}+\mathbf{B}.\end{equation}

Let $m=A_1$, $p=B_1$, and  let $(c_1+1,\ldots,c_m+1)$ and $(r_1+1,\ldots,r_p+1)$ be partitions defined by
$$(c_1+1,\ldots,c_m+1):=\overline{\mathbf{A}},$$
$$(r_1+1,\ldots,r_p+1):=\overline{\mathbf{B}}.$$
Let denote by $\mathbf{g}^i:=\overline{\mathbf{d}^i}$ with $\mathbf{g}^i=(g^i_1,g^i_2,\ldots)$, and $\mathbf{a}^i:=\overline{\mathbf{t}^i}$ with $\mathbf{a}^i=(a^i_1,a^i_2,\ldots)$, $i=1,\ldots,k$. Let $n=\max\{t^i_1| i=1,\ldots,k\}$. Let $\psi_1,\ldots,\psi_k$ be distinct irreducible  homogeneous polynomials from $\F[\lambda,\mu]$, and let $\tilde{\alpha}:\tilde{\alpha}_1|\cdots|\tilde{\alpha}_n$ and $\tilde{\gamma}:\tilde{\gamma}_1|\cdots|\tilde{\gamma}_{n+m+p}$ be polynomial chains defined by  
 $$\tilde{\alpha}_i=\psi_1^{a^1_{n+1-i}}\psi_2^{a^2_{n+1-i}}\ldots\psi_k^{a^k_{n+1-i}},\quad i=1,\ldots,n,$$
 and
 $$\tilde{\gamma}_i=\psi_1^{g^1_{n+m+p+1-i}}\psi_2^{g^2_{n+m+p+1-i}}\ldots\psi_k^{g^k_{n+m+p+1-i}},\quad i=1,\ldots,n+m+p.$$

Condition (\ref{?1}) is equivalent to \begin{equation}\label{q1} g^i_j\ge a^i_j,\quad i=1,\dots,k,\quad j\ge 1.\end{equation}
Also, (\ref{*1}) implies that for all $i=1,\ldots,k,$ we have $\max\{d^i_j-t^i_j|j\ge 1\}\le A_1+B_1=m+p$. Thus, 
 $\max\{\overline{g^i}_j-\overline{a^i}_j|j\ge 1\}\le m+p$, i.e. 
 \begin{equation}a^i_j\ge g^{i}_{j+m+p},\quad i=1,\ldots,k,\quad j\ge 1.\label{qq1}\end{equation} 
 
 Hence, (\ref{q1}) and (\ref{qq1}) give
 \begin{equation}\tilde{\gamma}_i|\tilde{\alpha}_i|\tilde{\gamma}_{i+m+p},\quad i=1,\ldots,n.\label{oo}\end{equation}

 Then by (\ref{1u})  we have
 \begin{equation}\label{1uy}(d(\tilde{\sigma}_{m+p}(\tilde{\alpha},\tilde{\gamma})),\ldots, d(\tilde{\sigma}_{1}(\tilde{\alpha},\tilde{\gamma})))=\overline{\overline{\mathbf{g}^1}-\overline{\mathbf{a}^1}}+\cdots+\overline{\overline{\mathbf{g}^k}-\overline{\mathbf{a}^k}}.\end{equation}
 Therefore  (\ref{*1}) becomes
\begin{equation} (c_1+1,\ldots,c_m+1)\cup(r_1+1,\ldots,r_p+1)\prec(d(\tilde{\sigma}_{m+p}(\tilde{\alpha},\tilde{\gamma})),\ldots, d(\tilde{\sigma}_{1}(\tilde{\alpha},\tilde{\gamma}))).\label{oo1}\end{equation}

  
  So by Theorem \ref{T}, since (\ref{oo}) and (\ref{oo1}) are valid, we have that there exists a polynomial chain $\tilde{\beta}:\tilde{\beta}_1|\cdots|\tilde{\beta}_{n+m}$ satisfying (1)--(4). Since $\tilde{\beta}_{n+m}|\tilde{\gamma}_{n+m+p}$, the only irreducible factors of $\tilde{\beta}_i$'s are $\psi_1,\ldots,\psi_k$. Let $\mathbf{b}^i=(b^i_1,\ldots,b^i_{n+m})$ be the corresponding partitions of $\psi_i$-elementary divisors of $\tilde{\beta}$, $i=1,\ldots,k$, i.e. 
   $$\tilde{\beta}_i=\psi_1^{b^1_{n+m+1-i}}\psi_2^{b^2_{n+m+1-i}}\ldots\psi_k^{b^k_{n+m+1-i}},\quad i=1,\ldots,n+m.$$
 Then  (1)--(4) imply:
\begin{eqnarray*}&{g}^i_j\ge {b}^i_j\ge {a}^i_j,\quad i=1,\ldots,k,\quad j\ge 1,\\
&(c_1+1,\ldots,c_m+1)\prec\overline{\overline{\mathbf{b}^1}-\overline{\mathbf{a}^1}}+\cdots+\overline{\overline{\mathbf{b}^k}-\overline{\mathbf{a}^k}}\\
&(r_1+1,\ldots,r_p+1)\prec\overline{\overline{\mathbf{g}^1}-\overline{\mathbf{b}^1}}+\cdots+\overline{\overline{\mathbf{g}^k}-\overline{\mathbf{b}^k}}.\end{eqnarray*}
Let $\mathbf{f}^i=\overline{\mathbf{b}^i}$, $i=1,\ldots,k$. Then duals of the conditions from above give
\begin{eqnarray*}&{d}^i_j\ge {f}^i_j\ge {t}^i_j,\quad i=1,\ldots,k,\quad j\ge 1,\\
&(\mathbf{f}^1-\mathbf{t}^1)\cup(\mathbf{f}^2-\mathbf{t}^2)\cup\cdots\cup(\mathbf{f}^k-\mathbf{t}^k)\prec \mathbf{A}\\
&(\mathbf{d}^1-\mathbf{f}^1)\cup(\mathbf{d}^2-\mathbf{f}^2)\cup\cdots\cup(\mathbf{d}^k-\mathbf{f}^k)\prec \mathbf{B},
\end{eqnarray*}
which proves Lemma \ref{nova}, as wanted. 
\kraj
\\
\begin{remark}\label{r}Since for if $k=1$ Lemma \ref{nova} reduces to Lemma \ref{silvaL}, we have that Theorem \ref{T} for $k=1$, i.e. in the case when $\gamma_{n+m+p}$ has only one irreducible factor, is  equivalent to Lemma \ref{silvaL}.\end{remark}

\vspace{0.5cm}

Since Theorem \ref{main} proves the equivalence between Lemma \ref{nova} and Theorem \ref{T}, as a corollary of Theorem \ref{main} we obtain that  Lemma \ref{nova} holds. It is a novel combinatorial result, that generalises Lemma \ref{silvaL}, whose applications in the control theory of linear systems are expected, and will be pursued in a future work.\\


\begin{example} It is well known that Theorem \ref{T} works over algebraically closed field (see e.g. \cite{cabral}). Let us comment how this translates into Lemmas \ref{silvaL} and \ref{nova}. 

The difference between arbitrary and algebraically closed field appears in the difference between (\ref{1}) and (\ref{1u}). Even in the case $k=1$ and $d(\psi_1)=2$, we would have
\begin{equation}\label{1uu}(d(\tilde{\sigma}_{m+p}(\tilde{\alpha},\tilde{\gamma})),\ldots, d(\tilde{\sigma}_{1}(\tilde{\alpha},\tilde{\gamma})))=2\,\,\overline{\overline{\mathbf{g}^1}-\overline{\mathbf{a}^1}}.\end{equation}
The analog of Lemma \ref{silvaL} that would be required in this case would be:\\

If $\mathbf{d}$, $\mathbf{t}$, $\mathbf{A}$ and $\mathbf{B}$ are partitions with $d_i\ge t_i$, $i\ge 1$, such that
\begin{equation}\label{gl}2(\mathbf{d}-\mathbf{t})\prec \mathbf{A+B}\end{equation} then there exists a partition $\mathbf{f}$
such that \begin{eqnarray*}
&d_i\ge f_i\ge t_i, \quad i\ge 1,\\
&2(\mathbf{f}-\mathbf{t})\prec \mathbf{A}\\
&2(\mathbf{d}-\mathbf{f})\prec \mathbf{B}.\end{eqnarray*}
However, this can be easily seen to be false. For example, let 
$\mathbf{A}=(1,1)$, $\mathbf{B}=(1,1)$, $d_1=t_1+1$ and $d_2=t_2+1$. 
Then (\ref{gl}) is satisfied, while there is no
$\mathbf{f}=(f_1,f_2)$ such that 
\begin{eqnarray}
&&d_i\ge f_i\ge t_i, \quad i\ge 1,\\
&&(2(f_1-t_1),2(f_2-t_2))\prec (1,1)\\
&&(2(d_1-f_1),2(d_2-f_2))\prec (1,1).\end{eqnarray}
\end{example}


\footnotesize
\begin{thebibliography}{20}


\bibitem{cabral} I. Cabral, F. C. Silva, Similarity invariants of completions of submatrices, Linear Algebra Appl. 169 (1992) 151-161.

\bibitem{super} M. Dodig, M. Sto\v si\'c, Similarity class of a matrix with prescribed submatrix, Linear and Multilinear Algebra, 57 (2009) 217-245.
\bibitem{pol} M. Dodig, M. Sto\v si\'c, Combinatorics of polynomial paths,  submitted. 


 \bibitem{silvaL} M. Dodig, F. C. Silva, Controllability of series connections of arbitrarily many linear systems, Linear Algebra Appl. 429 (2008) 122-141.

\bibitem{rc} M. Dodig, Completion of quasi-regular matrix pencils, Linear Algebra and its Applications 501 (2016) 198-241.

\bibitem{woat} M. Dodig, Descriptor systems under feedback and output injection, Operator Theory, Operator Algebras, and Matrix Theory, 2018.


















\bibitem{face} I. Gohberg, M. A. Kaashoek, F. van Schagen, Eigenvalues of completions of submatrices, Linear and Multilinear Algebra, 25 (1989)
55-70.



\bibitem{ped} G. Hardy, J. E. Littlewood, G. P\'olya, Inequalities, Cambridge University Press, 1991.

 










\bibitem{oliveira}
G.\ N.\ de Oliveira, Matrices with prescribed characteristic polynomial and a prescribed submatrix III, Monatsh. Math. 75 (1971), 441-446.


\bibitem{sa} E. M. S\'a, Imbedding conditions for $\lambda$ -matrices, Linear Algebra Appl. 24 (1979) 33-50.

\bibitem{tomp} R. C. Thompson, Interlacing inequalities for invariant factors,  Linear Algebra Appl. 24 (1979) 1-31.


\bibitem{zabala} I.\ Zaballa, Matrices with prescribed rows and invariant factors, Linear Algebra Appl.
87 (1987) 113-146.


\end{thebibliography}
\end{document}